\documentstyle{article}
\begin{document}
\input{amssym}
\author{Mehdi Nadjafikhah\thanks{Faculty of Pure Mathematics, Department of Mathematics, Iran University
of Science and technology, I.R. IRAN. e-mail:
m\_nadjafikhah@iust.ac.ir} \and Ahmad-Reza
Forough\thanks{a\_forough@iust.ac.ir}}
\title{Galilean geometry of motions}
\date{}
\maketitle
\abstract{In this paper we show that Galilean group is a matrix
Lie group and find its structure. Then provide the invariants of
special Galilean geometry of motions, by Olver's method of moving
coframes, we also find the corresponding $\{e\}-$structure.


\medskip \noindent A.M.S. 2000 Subject Classification: 58D19, 70E15.}
\section{Introduction:}
The method of moving coframes is one of the cornerstones of our
investigations for finding the invariants of the equivalence
problem which is done in \cite{F-O-1} and \cite{F-O-2}. Let us
summarize the basic method of moving coframes.

The basic steps are:
\begin{itemize}
\item[(i)] Determine the moving frame of order zero, by choosing a base
point and solving for Galilean group action.
\item[(ii)] Determine the invariant forms in this case (the finite
dimensional), they are the Maurer-Cartan forms, which computed by
direct use of the transformation group formulae but not the matrix
approach.
\item[(iii)] Use the invariant lift to pull-back the invariant forms,
leading to the moving coframe of order zero.
\item[(iv)] Determine the lifted invariants by finding the linear
dependencies among the restricted to horizontal components of the
moving coframe forms.
\item[(v)] Normalize any group-dependent invariants to convenient
constant values by solving for some of the unspecified parameters.
\item[(vi)] Successively eliminate parameters by substituting the
normalization formulae into the moving coframe and recomputing
dependencies.
\item[(vii)] After the parameters have all been normalized, the
differential invariants will appear through any remaining
dependencies among the final moving coframe elements. The
Invariant differential operators are found as the dual
differential operators to a basis for the invariant coframe forms.
\end{itemize}
At first we define some basic prerequisite from Galilean group.
Explanatory details are found in  \cite{L-B} , \cite{L}.
\paragraph{Definition 1.} The {Galilean group} is defined as
\begin{eqnarray*}
{\rm Gal}(3) &=& \left\{
\left[\begin{array}{ccc} 1& 0 & s \\
{\bf v} & R & {\bf y} \\ 0 & 0 & 1
\end{array}\right] \Bigg|\; s\in {\Bbb R}\,,\,{\bf y},{\bf v}\in{\Bbb R}^3\,,\,R\in {\rm O}(3)\right\}
\end{eqnarray*}
with a natural closed Lie subgroup structure of ${\rm GL}(5;{\Bbb
R})$ as
\begin{eqnarray*}
&& \hspace{-3.5cm} \left[ \begin{array}{ccc} 1&0&s_{{1}}\\{\bf
v}_{{1}}&R_{{1}}&{\bf y}_{{1}}\\0&0&1\end {array} \right]\,.\,
\left[
\begin {array}{ccc} 1&0&s_{{2}}\\{\bf
v}_{{2}}&R_{ {2}}&{\bf y}_{{2}}\\0&0&1\end {array} \right] = \nonumber \\
&=& \left[
\begin {array}{ccc}
1&0&s_1+s_2\\{\bf v}_{{1}}+R_{{1}}{\bf
v}_{{2}}&R_{{1}}R_{{2}}&{\bf y}_1+s_2{\bf v}_{{1}}+R_{1}{\bf
y}_{2} \\
0&0&1\end {array} \right] \\
 \left[\begin{array}{ccc} 1& 0 & s \\
{\bf v} & R & {\bf y} \\ 0 & 0 & 1
\end{array}\right]^{-1} &=& \left[ \begin {array}{ccc} 1  & 0 & -s \\
-R^{-1}{\bf v} & R^{-1} & R^{-1}(s{\bf v}-{\bf y}) \\
0 & 0 & 1
\end {array} \right]. \nonumber
\end{eqnarray*}
This is a $10-$dimensional Lie group. The {\it special Galilean
group} is defined as connected component of $e$ in ${\rm Gal}(3)$
and denoted by ${\rm SGal}(3)$.
\paragraph{Definition 2.} Let we identify the ${\Bbb R}^4$ by
\begin{eqnarray*}
{\Bbb R}^4 &=& \left\{
\left[\begin{array}{c} t \\
{\bf x} \\ 1
\end{array}\right] \Bigg|\; t\in {\Bbb R}\,,\,{\bf x}=\left[\begin{array}{c} x \\
y \\ z \end{array}\right] \in{\Bbb R}^3\right\}
\end{eqnarray*}
with the natural $4-$manifold structure. Then, we can define
naturally the smooth action of ${\rm Gal}(3)$ on ${\Bbb R}^4$ as
\begin{eqnarray*}
\left[\begin{array}{ccc} 1& 0 & s \\
{\bf v} & R & {\bf y} \\ 0 & 0 & 1
\end{array}\right] \bullet
\left[\begin{array}{c} t \\
{\bf x} \\ 1
\end{array}\right] = \left[\begin{array}{c} t+s \\ R{\bf x}+t{\bf v}+{\bf y} \\ 1
\end{array}\right]
\end{eqnarray*}

By elementary algebraic computations, we find the structure of
special Galilean group,
\paragraph{Theorem 1.}
{\it Let
\begin{eqnarray*}
G_1 &=& \left\{
\left[\begin{array}{ccc} 1& 0 & 0 \\
{\bf v} & R & {\bf 0} \\ 0 & 0 & 1
\end{array}\right] \Bigg|\; {\bf v}\in{\Bbb R}^3\,,\,R\in {\rm
O}(3)\right\} \nonumber \leq {\rm GL}(4;{\Bbb R})
\end{eqnarray*}
be the group of uniformly special Galilean  motions,
\begin{eqnarray*}
G_2 &=& \left\{
\left[\begin{array}{ccc} 1& 0 & s \\
{\bf 0} & I_3 & {\bf y} \\ 0 & 0 & 1
\end{array}\right] \Bigg|\; s\in{\Bbb R}\,,\,{\bf
y}\in {\Bbb R}^3 \right\} \cong ({\Bbb R}^4,+)
\end{eqnarray*}
be the group of shifts of origin,
\begin{eqnarray*}
G_3 &=& \left\{
\left[\begin{array}{ccc} 1& 0 & 0 \\
{\bf 0} & R & {\bf 0} \\ 0 & 0 & 1
\end{array}\right] \Bigg|\; R\in {\rm
SO}(3)\right\} \cong {\rm SO}(3)
\end{eqnarray*}
be the group of rotations of reference frame, and
\begin{eqnarray*}
G_4 &=& \left\{
\left[\begin{array}{ccc} 1& 0 & 0 \\
{\bf v} & I_3 & {\bf 0} \\ 0 & 0 & 1
\end{array}\right] \Bigg|\; {\bf
v}\in {\Bbb R}^3 \right\} \cong ({\Bbb R}^3,+)
\end{eqnarray*}
be the group of uniformly frame motions. Then, $G_2 \unlhd {\rm
SGal}(3)$, ${\rm SGal}(3) \cong G_1\propto G_2$,
 $G_4 \unlhd G_1$, $G_1 \cong G_3\propto G_4$, and
 ${\rm SGal}(3) \cong ({\rm SO}(3)\propto {\Bbb R}^3)\propto {\Bbb
R}^4$.}

\medskip In the following theorem, we explain the algebraic
structure of the infinitesimal group action
$\widehat{\goth{Gal}}(3)$ induced by the action ${\rm SGal}(3)$ on
${\Bbb R}^4$,
\paragraph{Theorem 2.}
{\it The Lie algebra of infinitesimal group action
$\widehat{\goth{Gal}}(3)={\rm Span}_{\Bbb
R}\{\widehat{X}_1,\cdots,\widehat{X}_{10}\}$  induced by the
action ${\rm SGal}(3)$ on ${\Bbb R}^4$,  has infinitesimal
generators:
\begin{eqnarray*}
\widehat{X}_1=\partial_t,\;\;
\begin{array}{l} \widehat{X}_2=\partial_x, \\
\widehat{X}_3=\partial_y, \\ \widehat{X}_4=\partial_z,
\end{array}\;\;
\begin{array}{l} \widehat{X}_5=t\,\partial_x, \\
\widehat{X}_6=t\,\partial_y, \\ \widehat{X}_7=t\,\partial_z,
\end{array}\;\;
\begin{array}{l}
\widehat{X}_8=y\,\partial_x-x\,\partial_y, \\
\widehat{X}_9=x\,\partial_z-z\,\partial_x,  \\
\widehat{X}_{10}=z\,\partial_y-y\,\partial_z, \end{array}
\end{eqnarray*}
with the following structure:
\begin{eqnarray*}
\begin{array}{|c||c|ccc|ccc|ccc|}
\hline &&&&&&&&&& \\[-1mm]
      \!\!&\!\!\widehat X_1   \!\!&\!\!\widehat X_2  \!\!&\!\!\widehat X_3 \!\!&\!\!\widehat X_4 \!\!&\!\!\widehat X_5 \!\!&\!\!\widehat X_6 \!\!&\!\!\widehat X_7 \!\!&\!\!\widehat X_8    \!\!&\!\!\widehat X_9    \!\!&\!\!\widehat X_{10}
      \\[2mm]
\hline \hline &&&&&&&&&& \\[-1mm]
\widehat X_1   \!\!&\!\! 0    \!\!&\!\!  0  \!\!&\!\! 0  \!\!&\!\! 0  \!\!&\!\! \widehat X_2\!\!&\!\! \widehat X_3\!\!&\!\! \widehat X_4\!\!&\!\! 0     \!\!&\!\!  0    \!\!&\!\!  0    \\[2mm]
\hline &&&&&&&&&& \\[-1mm]
\widehat X_2   \!\!&\!\! 0    \!\!&\!\!  0  \!\!&\!\! 0  \!\!&\!\! 0  \!\!&\!\! 0  \!\!&\!\! 0  \!\!&\!\! 0  \!\!&\!\! -\widehat X_3  \!\!&\!\! \widehat X_4   \!\!&\!\!  0    \\[2mm]
\widehat X_3   \!\!&\!\! 0    \!\!&\!\!  0  \!\!&\!\! 0  \!\!&\!\! 0  \!\!&\!\! 0  \!\!&\!\! 0  \!\!&\!\! 0  \!\!&\!\! \widehat X_2   \!\!&\!\! 0     \!\!&\!\! -\widehat X_4  \\[2mm]
\widehat X_4   \!\!&\!\! 0    \!\!&\!\!  0  \!\!&\!\! 0  \!\!&\!\! 0  \!\!&\!\! 0  \!\!&\!\! 0  \!\!&\!\! 0  \!\!&\!\! 0     \!\!&\!\!-\widehat X_2   \!\!&\!\!  \widehat X_3  \\[2mm]
\hline &&&&&&&&&& \\[-1mm]
\widehat X_5   \!\!&\!\! -\widehat X_2 \!\!&\!\!  0  \!\!&\!\! 0  \!\!&\!\! 0  \!\!&\!\! 0  \!\!&\!\! 0  \!\!&\!\! 0  \!\!&\!\! -\widehat X_6  \!\!&\!\! \widehat X_7   \!\!&\!\!  0    \\[2mm]
\widehat X_6   \!\!&\!\! -\widehat X_3 \!\!&\!\!  0  \!\!&\!\! 0  \!\!&\!\! 0  \!\!&\!\! 0  \!\!&\!\! 0  \!\!&\!\! 0  \!\!&\!\! \widehat X_5   \!\!&\!\! 0     \!\!&\!\! -\widehat X_7  \\[2mm]
\widehat X_7   \!\!&\!\! -\widehat X_4 \!\!&\!\!  0  \!\!&\!\! 0  \!\!&\!\! 0  \!\!&\!\! 0  \!\!&\!\! 0  \!\!&\!\! 0  \!\!&\!\! 0     \!\!&\!\!-\widehat X_5   \!\!&\!\! \widehat X_6   \\[2mm]
\hline &&&&&&&&&& \\[-1mm]
\widehat X_8   \!\!&\!\! 0    \!\!&\!\! \widehat X_3 \!\!&\!\!-\widehat X_2\!\!&\!\! 0  \!\!&\!\! \widehat X_6\!\!&\!\!-\widehat X_5\!\!&\!\! 0  \!\!&\!\! 0     \!\!&\!\!-\widehat X_{10}\!\!&\!\! \widehat X_9   \\[2mm]
\widehat X_9   \!\!&\!\! 0    \!\!&\!\!-\widehat X_4 \!\!&\!\! 0  \!\!&\!\! \widehat X_2\!\!&\!\!-\widehat X_7\!\!&\!\! 0  \!\!&\!\! \widehat X_5\!\!&\!\! \widehat X_{10}\!\!&\!\! 0     \!\!&\!\!-\widehat X_8   \\[2mm]
\widehat X_{10}\!\!&\!\! 0    \!\!&\!\! 0   \!\!&\!\!\widehat X_4 \!\!&\!\!-\widehat X_3\!\!&\!\! 0  \!\!&\!\! \widehat X_7\!\!&\!\!-\widehat X_6\!\!&\!\!-\widehat X_9   \!\!&\!\! \widehat X_8   \!\!&\!\! 0     \\[2mm]
\hline
\end{array}
\end{eqnarray*}}
\section{Computation of Maurer-Cartan forms}
\paragraph{Definition 3.}
By multiplying 3 rotations
\begin{eqnarray*}
\left[ \begin {array}{ccc} \cos \theta_3 &\sin \theta_3 &0\\-\sin
\theta_3 &\cos \theta_3 &0
\\0&0&1\end {array} \right],\;
\left[ \begin {array}{ccc} \cos \theta_2 &0&\sin \theta_2
\\0&1&0
\\-\sin \theta_2 &0&\cos
\theta_2 \end {array} \right],\; \left[ \begin {array}{ccc}
1&0&0\\0&\cos \theta_1 &\sin \theta_1
\\0&-\sin \theta_1 &\cos
\theta_1\end {array} \right]
\end{eqnarray*}
respectively about $x$, $y$ and $z-$axis, we define
\begin{eqnarray*}
R= \tiny \left[ \begin {array}{ccc} \cos \theta_2 \cos \theta_3
&\cos \theta_1 \sin
 \theta_3 -\sin \theta_1 \sin
 \theta_2 \cos \theta_3 &\sin
 \theta_1 \sin \theta_3 +\cos
 \theta_1 \sin \theta_2 \cos
 \theta_3 \\-\cos \theta_2 \sin \theta_3 &\cos \theta_1 \cos \theta_3 +\sin \theta_1 \sin \theta_2 \sin
\theta_3 &\sin \theta_1 \cos \theta_3 -\cos\theta_1 \sin \theta_2
\sin \theta_3
\\- \sin \theta_2 &-\sin
\theta_1 \cos \theta_2 &\cos \theta_1 \cos \theta_2
\end {array} \right]. \end{eqnarray*}
which is an arbitrary element of ${\rm SO}(3)$.

\medskip In order to determine the Maurer-Cartan forms, we would
rather use the direct method,  more details found in page 10 of
\cite{F-O-1}.
\paragraph{Theorem 3.}
{\it The independent Maurrer-Cartan $1-$forms of ${\rm SGal}(3)$
are
\begin{eqnarray*}
\mu_1 &=& ds,\\
\mu_2 &=& \sin\theta_2\,d\theta_1 - d\theta_3,\\
\mu_3 &=& \cos\theta_2\cos\theta_3\,d\theta_1 - \sin\theta_3\,d\theta_2,\\
\mu_4 &=& \cos\theta_2\sin\theta_3\,d\theta_1 + \cos\theta_3\,d\theta_2,\\
\mu_5 &=& \cos\theta_1\cos\theta_2\,dv_1 - \sin\theta_1\cos\theta_2\,dv_2 - \sin\theta_2\,dv_3,\\
\mu_6 &=& \big(\cos\theta_1\cos\theta_2v_1 - \sin\theta_1\cos\theta_2v_2 - \sin\theta_2v_3\big)\,ds\\
      &  & - \cos\theta_1\cos\theta_2\,dy_1 + \sin\theta_1\cos\theta_2\,dy_2,\\
\mu_7 &=& - \big(\sin\theta_1\sin\theta_3
           + \cos\theta_1\sin\theta_2\cos\theta_3\big)\,dv_1\\
      &  & + \big(\sin\theta_1\sin\theta_2\cos\theta_3 - \cos\theta_1\sin\theta_3\big)\,dv_2 - \cos\theta_2\cos\theta_3\,dv_3,\\
\mu_8 &=& \big(\cos\theta_1\sin\theta_2\sin\theta_3
           - \sin\theta_1\cos\theta_3\big)\,dv_1\\
      &  & - \big(\sin\theta_1\sin\theta_2\sin\theta_3 + \cos\theta_1\cos\theta_3\big)\,dv_2 + \cos\theta_2\sin\theta_3\,dv_3,\\
\mu_9 &=&
            \Big((\sin\theta_1\sin\theta_3+\cos\theta_1\sin\theta_2\cos\theta_3)\,v_1
            +(\cos\theta_1\sin\theta_3\\
      &  & -\sin\theta_1\sin\theta_2\cos\theta_3)v_2+\cos\theta_2\cos\theta_3v_3\Big)\,ds\\
      &  & -\big(\cos\theta_1\sin\theta_2\cos\theta_3+\sin\theta_1\sin\theta_3\big)\,dy_1\\
      &  & +\big(-\cos\theta_1\sin\theta_3+\sin\theta_1\sin\theta_2\cos\theta_3\big)\,dy_2 -
      \cos\theta_2\cos\theta_3\,dy_3,\\
\mu_{10}&=&
            \Big((\sin\theta_1\cos\theta_3-\cos\theta_1\sin\theta_2\sin\theta_3)\,v_1
            +(\cos\theta_1\cos\theta_3\\
       &  & -\sin\theta_1\sin\theta_2\sin\theta_3)v_2+\cos\theta_2\sin\theta_3v_3\Big)\,ds\\
       &  & -\big(\cos\theta_1\sin\theta_2\cos\theta_3+\sin\theta_1\sin\theta_3\big)\,dy_1\\
       &  & +\big(\sin\theta_1\sin\theta_2\cos\theta_3-\cos\theta_1\sin\theta_3\big)\,dy_2 -
       \cos\theta_2\cos\theta_3\,dy_3.
\end{eqnarray*}}

\medskip \noindent {\it Proof:} Given $g\in {\rm SGal}(3)$ and ${\bf z}\in {\Bbb R}^4$, we explicitly
write the group transformation $\bar{\bf z}=g\cdot {\bf z}$ in
coordinate form:
\begin{eqnarray*}
\bar{z}_1 &=& H^1({\bf z},g) = t+s\\
\bar{z}_2 &=& H^2({\bf z},g)\\
          &=& {v_1}\,t+(\cos  {\theta_2}  \cos  {\theta_1} ) {x_1} + ( \cos
{\theta_3}  \sin  {\theta_1}
  -\sin  {\theta_3}  \sin  {\theta_2}
  \cos  {\theta_1}   ) {x_2}\\
          & & + (
\sin  {\theta_3}  \sin  {\theta_1}  + \cos  {\theta_3}  \sin
{\theta_2}  \cos
  {\theta_1}   ) {x_3}+{y_1}\\
\bar{z}_3 &=& H^3({\bf z},g)\\
          &=&
{v_2}\,t-(\cos  {\theta_2}  \sin
  {\theta_1})  {x_1}+ ( \cos  {\theta_3}
  \cos  {\theta_1}  +\sin  {\theta_3}
  \sin  {\theta_2}  \sin  {\theta_1}
   ) {x_2}\\
          & & + ( \sin  {\theta_3}
\cos  {\theta_1}  -\cos  {\theta_3}  \sin  {\theta_2}  \sin
{\theta_1}
 ) {x_3}+{y_2}\\
\bar{z}_4 &=& H^4({\bf z},g)\\
          &=& {v_3}\,t-\sin
{\theta_2}  {x_1}-(\sin  {\theta_3}  \cos
  {\theta_2})  {x_2}+(\cos  {\theta_3}
  \cos  {\theta_2} ) {x_3}+{y_3}.
\end{eqnarray*}
We then compute the differentials of the group transformations:
$$ d\bar{z}_i = \sum_{k=1}^4 \frac{\partial H^i}{\partial z_k}\,dz_k
+ \sum_{j=1}^{10}\frac{\partial H^i}{\partial
g^j}\,dg^j,\hspace{1cm} i=1,\cdots,4, $$ or more compactly
\begin{eqnarray}
d\bar{\bf z}=H_{\bf z}\,d{\bf z}+H_g\,dg. \label{olver:1}
\end{eqnarray}
 Next, set $d\bar{\bf
z}=0$ in (\ref{olver:1}), and solve the resulting system of linear
equations for the differentials $dz_k$. This leads to the formulae
\begin{eqnarray*}
-d{\bf z} = F\,dg = (H_{\bf z}^{-1}\cdot H_g)\,dg,
\end{eqnarray*}
or, in full detail,
\begin{eqnarray}
-dz_k=\sum_{j=1}^{10} F_j^k({\bf z},g)\,dg^j,\hspace{1cm}
i=1,\cdots,4. \label{olver:2}
\end{eqnarray}
Then, for each $k$ and each fixed ${\bf z}_0\in {\Bbb R}^4$, the
one-form $\mu_0=\sum_{j=1}^{10}F_j^k({\bf z}_0,g)\,dg^j$ is a
left-invariant Maurer-Cartan form on the group ${\rm SGal}(3)$.
Alternatively, if one expands the right hand side of
(\ref{olver:2}) in power series in $\bf z$,
\begin{eqnarray*}
\sum_{j=1}^{10} F_j^k({\bf z},g)\,dg^j = \sum_{i=0}^\infty
z_i\,\mu_i,
\end{eqnarray*}
then each coefficient $\mu_i$ also forms a left-invariant
Maurer-Cartan form on ${\rm SGal}(3)$. \hfill\ $\Box$
\section{Zero order moving coframes}
  Throughout this paper, we remind you that  ${\rm SGal}(3)$ is a
$10-$dimensional Lie group and $M$ is a $4-$dimensional manifold.
\paragraph{Definition 4.}
A smooth map $\rho^{(0)} :M\to{\rm SGal}(3)$ is called a
compatible lift with base point $z_0$ if it satisfies
$$\rho^{(0)} (z).z_0 =z,\;\;\;\;\;\; z\in M.$$

\medskip Now, let $\rho^{(0)} :M\to {\rm SGal}(3)$ be a compatible
lift with base point ${\bf z}_0=[0,{\bf 0},1]^T\in{\Bbb R}^4$,
then we have, $s=t$ and ${\bf y}={\bf x}$. Thus,
\paragraph{Theorem 4.}
{\it The most general zero order compatible lift has the form
\begin{eqnarray*}
\rho^{(0)}(t,{\bf x};{\bf v},{\bf \theta}) = \left[\begin{array}{ccc} 1& 0 & t \\
{\bf v} & R & {\bf x} \\ 0 & 0 & 1
\end{array}\right].
\end{eqnarray*}}

The next step is to characterize the group transformations by a
collection of differential forms. In the finite-dimensional
situation that we are currently considering, these will be
obtained by pulling back the left-invariant Maurer-Cartan forms
$\bf \mu$ on ${\rm SGal}(3)$ to the order zero moving frame bundle
${\cal B}_0$ using the compatible lift.

The resulting one-forms ${\bf \zeta}^{(0)}=\rho^{(0)*}{\bf \mu}$
will provide an invariant coframe on ${\cal B}_0$, which we name
the {\it moving coframe of the zero order}. The moving coframe
forms ${\bf \zeta}^{(0)}$ clearly satisfy the same Maurer-Cartan
structure equations. Thus
\paragraph{Theorem 5.} {\it The zero order moving coframe is
\begin{eqnarray*}
\zeta^{(0)}_1 &=& dt,\\
\zeta^{(0)}_2 &=& \sin\theta_2\,d\theta_1 - d\theta_3,\\
\zeta^{(0)}_3 &=& \cos\theta_2\cos\theta_3\,d\theta_1 - \sin\theta_3\,d\theta_2,\\
\zeta^{(0)}_4 &=& \cos\theta_2\sin\theta_3\,d\theta_1 + \cos\theta_3\,d\theta_2,\\
\zeta^{(0)}_5 &=& \cos\theta_1\cos\theta_2\,dv_1 - \sin\theta_1\cos\theta_2\,dv_2 - \sin\theta_2\,dv_3,\\
\zeta^{(0)}_6 &=& \big(\cos\theta_1\cos\theta_2v_1 - \sin\theta_1\cos\theta_2v_2 - \sin\theta_2v_3)\big)\,dt\\
      &  & - \cos\theta_1\cos\theta_2\,dx_1 + \sin\theta_1\cos\theta_2\,dx_2,\\
\zeta^{(0)}_7 &=& - \big(\sin\theta_1\sin\theta_3
           + \cos\theta_1\sin\theta_2\cos\theta_3\big)\,dv_1\\
      &  & + \big(\sin\theta_1\sin\theta_2\cos\theta_3 - \cos\theta_1\sin\theta_3\big)\,dv_2 - \cos\theta_2\cos\theta_3\,dv_3\\
\zeta^{(0)}_8 &=& \big(\cos\theta_1\sin\theta_2\sin\theta_3
           - \sin\theta_1\cos\theta_3\big)\,dv_1,\\
      &  & - \big(\sin\theta_1\sin\theta_2\sin\theta_3 + \cos\theta_1\cos\theta_3\big)\,dv_2 + \cos\theta_2\sin\theta_3\,dv_3,\\
\zeta^{(0)}_9 &=&
            \Big((\sin\theta_1\sin\theta_3+\cos\theta_1\sin\theta_2\cos\theta_3)\,v_1
            +(\cos\theta_1\sin\theta_3\\
      &  & -\sin\theta_1\sin\theta_2\cos\theta_3)v_2+\cos\theta_2\cos\theta_3v_3\Big)\,dt\\
      &  & -\big(\cos\theta_1\sin\theta_2\cos\theta_3+\sin\theta_1\sin\theta_3\big)\,dx_1\\
      &  & +\big(-\cos\theta_1\sin\theta_3+\sin\theta_1\sin\theta_2\cos\theta_3\big)\,dx_2 -
      \cos\theta_2\cos\theta_3\,dx_3,\\
\zeta^{(0)}_{10}&=&
            \Big((\sin\theta_1\cos\theta_3-\cos\theta_1\sin\theta_2\sin\theta_3)\,v_1
            +(\cos\theta_1\cos\theta_3\\
       &  & -\sin\theta_1\sin\theta_2\sin\theta_3)v_2+\cos\theta_2\sin\theta_3v_3\Big)\,dt\\
       &  & -\big(\cos\theta_1\sin\theta_2\cos\theta_3+\sin\theta_1\sin\theta_3\big)\,dx_1\\
       &  & +\big(\sin\theta_1\sin\theta_2\cos\theta_3-\cos\theta_1\sin\theta_3\big)\,dx_2 -
       \cos\theta_2\cos\theta_3\,dx_3,
\end{eqnarray*}
which forms a basis for the space of one-forms on ${\cal
B}_0={\Bbb R}^4\times G_1$ \hfill\ $\Box$.}
\section{First order moving coframes}
\paragraph{Definition 5.}
A motion is a curve coincides with the graph of a function ${\bf
x} = {\bf x}(t):{\Bbb R}\to{\Bbb R}^3$.

\medskip We restrict the moving coframe forms to the motion (curve),
which amounts to replacing the differential $d{\rm x}$  by its
"horizontal" component ${\rm x}_t\,dt$. If we interpret the
derivative ${\rm x}_t$ as a coordinate on the first jet space $J^1
= J^1({\Bbb R}^1;{\Bbb R}^3)\cong{\Bbb R}^7$ of motions in ${\Bbb
R}^4$, then the restriction of a differential form to the motion
can be reinterpreted as the natural projection of the one-form
$d{\rm x}$ on $J^1$ to its horizontal component, using the
canonical decomposition of differential forms on the jet space
into horizontal and contact components. Indeed, the vertical
component of the form $d{\rm x}$ is the contact form $d{\bf
x}-{\bf x}_t\,dt$, which vanishes on all prolonged sections of the
first jet bundle $J^1({\Bbb R}^1;{\Bbb R}^3)$. Therefore,
\paragraph{Theorem 6.}
{\it The restricted (or horizontal) moving coframe forms are
defined on $7-$dimensional manifold $\{J^1{\rm x}=(t,{\rm
x}(t),{\rm x}_t(t))\}\times G_1\subset J^1{\cal B}_0$ and
explicitly given by $\eta^{(0)}_i = \zeta^{(0)}_i$, for
$i=1,2,3,4,5,7,8$, and their linear dependencies are
$\eta^{(0)}_6=j_1\,\eta^{(0)}_1$, $\eta^{(0)}_9=j_2\,\eta^{(0)}_1$
and $\eta^{(0)}_{10}=j_3\,\eta^{(0)}_1$, where
\begin{eqnarray*}
J_1 &=& -\cos\theta_1\cos\theta_2v_1
           +\sin\theta_1\cos\theta_2\,v_2 +\sin\theta_2v_3 \\
       & & +\cos\theta_1\cos\theta_2x_1'
       +\sin\theta_1\cos\theta_2x_2'-\sin\theta_2x_3',\\
J_2 &=&
\big(\cos\theta_1\sin\theta_2\sin\theta_3-\sin\theta_1\cos\theta_3\big)v_1\\
& &
-\big(\cos\theta_1\cos\theta_3+\sin\theta_1\sin\theta_2\sin\theta_3\big)v_2
+ \cos\theta_2\sin\theta_3v_3\\
&&+\big(\sin\theta_1\cos\theta_3
-\cos\theta_1\sin\theta_2\sin\theta_3 \big)x_1'\\
& & +\big(\sin\theta_1\sin\theta_2\sin\theta_3
+\cos\theta_1\cos\theta_3\big)x_2' -\cos\theta_2\sin\theta_3x_3',\\
J_3 &=&
-\big(\sin\theta_1\sin\theta_3+\cos\theta_1\sin\theta_2\cos\theta_3\big)v_1\\
& &
+\big(\sin\theta_1\sin\theta_2\cos\theta_3-\cos\theta_1\sin\theta_3\big)v_2
-\cos\theta_2\cos\theta_3v_3\\
& &
+\big(\sin\theta_1\sin\theta_3+\cos\theta_1\sin\theta_2\cos\theta_3\big)x_1'\\
& &
+\big(\cos\theta_1\sin\theta_3-\sin\theta_1\sin\theta_2\cos\theta_3\big)x_2'
+\cos\theta_2\cos\theta_3x_3'.
\end{eqnarray*}}

\medskip By assumptions $J_1=J_2=J_3=0$, we have ${\bf v}={\bf
x}_t$. Thus,
\paragraph{Theorem 7.}
{\it The first order compatible lift has the form:
\begin{eqnarray*}
\rho^{(1)}(t,{\bf x};{\bf x}_t,{\bf \theta}) = \left[\begin{array}{ccc} 1& 0 & t \\
{\bf x}_t & R & {\bf x} \\ 0 & 0 & 1
\end{array}\right].
\end{eqnarray*}}

The resulting one-forms ${\bf \zeta}^{(1)}=\rho^{(1)*}{\bf \mu}$
will provide an invariant coframe on ${\cal B}_1$, which we name
the {\it moving coframe of the first order}.
By substituting the map $\rho^{(1)}$ in $\zeta^{(0)}$ and
restricting to the first prolongation or jet of the motion, namely
${\bf x}={\bf x}(t)\,\,,{\bf x}_t={\bf x}'(t)$ we have,
\paragraph{Theorem 8.} {\it The first order moving coframe is
\begin{eqnarray*}
\zeta^{(1)}_1 &=& dt,\\
\zeta^{(1)}_2 &=& \sin\theta_2\,d\theta_1 - d\theta_3,\\
\zeta^{(1)}_3 &=& \cos\theta_2\cos\theta_3\,d\theta_1 - \sin\theta_3\,d\theta_2,\\
\zeta^{(1)}_4 &=& \cos\theta_2\sin\theta_3\,d\theta_1 + \cos\theta_3\,d\theta_2,\\
\zeta^{(1)}_5 &=& \cos\theta_1\cos\theta_2\,dx'_1 - \sin\theta_1\cos\theta_2\,dx'_2 - \sin\theta_2\,dx'_3,\\
\zeta^{(1)}_6 &=& \big(\cos\theta_1\cos\theta_2x'_1 - \sin\theta_1\cos\theta_2x'_2 - \sin\theta_2x'_3)\big)\,dt\\
      &  & - \cos\theta_1\cos\theta_2\,dx_1 + \sin\theta_1\cos\theta_2\,dx_2,\\
\zeta^{(1)}_7 &=& - \big(\sin\theta_1\sin\theta_3
           + \cos\theta_1\sin\theta_2\cos\theta_3\big)\,dx'_1\\
      &  & + \big(\sin\theta_1\sin\theta_2\cos\theta_3 - \cos\theta_1\sin\theta_3\big)\,dx'_2 - \cos\theta_2\cos\theta_3\,dx'_3\\
\zeta^{(1)}_8 &=& \big(\cos\theta_1\sin\theta_2\sin\theta_3
           - \sin\theta_1\cos\theta_3\big)\,dx'_1,\\
      &  & - \big(\sin\theta_1\sin\theta_2\sin\theta_3 + \cos\theta_1\cos\theta_3\big)\,dx'_2 + \cos\theta_2\sin\theta_3\,dx'_3,\\
\zeta^{(1)}_9 &=&
            \Big((\sin\theta_1\sin\theta_3+\cos\theta_1\sin\theta_2\cos\theta_3)\,x'_1
            +(\cos\theta_1\sin\theta_3\\
      &  & -\sin\theta_1\sin\theta_2\cos\theta_3)x'_2+\cos\theta_2\cos\theta_3x'_3\Big)\,dt\\
      &  & -\big(\cos\theta_1\sin\theta_2\cos\theta_3+\sin\theta_1\sin\theta_3\big)\,dx_1\\
      &  & +\big(-\cos\theta_1\sin\theta_3+\sin\theta_1\sin\theta_2\cos\theta_3\big)\,dx_2 -
      \cos\theta_2\cos\theta_3\,dx_3,\\
\zeta^{(1)}_{10}&=&
            \Big((\sin\theta_1\cos\theta_3-\cos\theta_1\sin\theta_2\sin\theta_3)\,x'_1
            +(\cos\theta_1\cos\theta_3\\
       &  & -\sin\theta_1\sin\theta_2\sin\theta_3)x'_2+\cos\theta_2\sin\theta_3x'_3\Big)\,dt\\
       &  & -\big(\cos\theta_1\sin\theta_2\cos\theta_3+\sin\theta_1\sin\theta_3\big)\,dx_1\\
       &  & +\big(\sin\theta_1\sin\theta_2\cos\theta_3-\cos\theta_1\sin\theta_3\big)\,dx_2 -
       \cos\theta_2\cos\theta_3\,dx_3,
\end{eqnarray*}
which is an invariant coframe on ${\cal B}_1=\{(t,{\rm x},{\rm
x}_t,\theta)\}\cong J^1({\Bbb R};{\Bbb R}^3)\times G_3\cong {\Bbb
R}^7\times {\rm SO}(3)$ \hfill\ $\Box$. }
\section{Second order moving coframes}
By restricting $\zeta^{(1)}$ to the second prolongation $J^2{\rm
x}\times G_3$, which is a four dimensional manifold, we have
\paragraph{Theorem 9.}
{\it The restricted (or horizontal) moving coframe forms are
explicitly given by
\begin{eqnarray*}
\eta^{(1)}_1 &=& dt,\\
\eta^{(1)}_2 &=& \sin\theta_2\,d\theta_1 - d\theta_3,\\
\eta^{(1)}_3 &=& \cos\theta_2\cos\theta_3\,d\theta_1 - \sin\theta_3\,d\theta_2,\\
\eta^{(1)}_4 &=& \cos\theta_2\sin\theta_3\,d\theta_1 +
\cos\theta_3\,d\theta_2,
\end{eqnarray*}
and their linear dependencies are, $\eta^{(1)}_5=J_1\eta^{(1)}_1$,
$\eta^{(1)}_7=J_2\eta^{(1)}_1$, $\eta^{(1)}_8=J_3\eta^{(1)}_1$,
and $\eta^{(1)}_6=\eta^{(1)}_9=\eta^{(1)}_{10}=0$, where
\begin{eqnarray*}
J_1 &=& (\cos\theta_1\sin\theta_2\sin\theta_3-\sin\theta_1\cos\theta_3)x_1''\\
    & &
    -(\sin\theta_1\sin\theta_2\sin\theta_3+\cos\theta_1\cos\theta_3)x_2''+\cos\theta_2\sin\theta_3x_3'',\\
J_2 &=& -(\sin\theta_1\sin\theta_3+\cos\theta_3\sin\theta_2\cos\theta_1)x_1''\\
    & & +(\sin\theta_1\sin\theta_2\cos\theta_3-\cos\theta_1\sin\theta_3)x_2''-\cos\theta_2\cos\theta_3x_3'',\\
J_3 &=&
       -\cos\theta_2\cos\theta_1x_1''+\cos\theta_2\sin\theta_1x_2''+\sin\theta_2x_3''.\\
\end{eqnarray*}}

If we assume $(J_1,J_2,J_3)=(-a,0,0)$, where the length of
acceleration $\|{\bf x}_{tt}\|$ is denoted by $a$ , then we have
$R{\bf x}_{tt}=(a,0,0)$, and by simple computations, have
\begin{eqnarray*}
\theta_1=-\arctan\Big(\frac{x_2''}{x_1''}\Big) \hspace{1cm}
\mbox{and} \hspace{1cm} \theta_2=\arcsin\Big(\frac{x_3''}{\|{\bf
x}_{tt}\|}\Big).
\end{eqnarray*}

It can be also easily seen that $a$ is an invariant.

Now we choose a cross section $K=\{t=0,{\rm x}={\rm 0},{\rm
x}_t={\rm 0},\|{\rm x}_{tt}\|=a,\theta={\rm 0}\}$. By recomputing
the forms ${\bf \zeta}^{(2)}=\rho^{(2)*}{\bf \mu}$, we have
\paragraph{Theorem 10.} {\it The second order moving coframe is
\begin{eqnarray*}
\zeta^{(2)}_1 \!\!\!\!&=&\!\!\!\! dt, \\
\zeta^{(2)}_2 \!\!\!\!&=&\!\!\!\!
-d\theta_3 + \frac{x''_2x''_3\,dx_1''}{a({x_1''}^2+{x_2''}^2)}-\frac{x''_1x''_3\,dx_2''}{a({x_1''}^2+{x_2''}^2)}, \\
\zeta^{(2)}_3 \!\!\!\!&=&\!\!\!\!
\frac{x_1''x_3''\sin\theta_3+ax_2''\cos\theta_3}{a^2\sqrt{{x_1''}^2+{x_2''}^2}}\,dx_1''
-
\frac{x_2''x_3''\sin\theta_3+ax_1''\cos\theta_3}{a^2\sqrt{{x_1''}^2+{x_2''}^2}}\,dx_2''
\\ & & +
\frac{\sin\theta_3\sqrt{{x_1''}^2+{x_2''}^2}}{a^2}\,dx_3'', \\
\zeta^{(2)}_4 \!\!\!\!&=&\!\!\!\!
-\frac{x_1''x_3''\cos\theta_3+ax_2''\sin\theta_3}{a^2\sqrt{{x_1''}^2+{x_2''}^2}}\,dx_1''
-
\frac{x_2''x_3''\cos\theta_3+ax_1''\sin\theta_3}{a^2\sqrt{{x_1''}^2+{x_2''}^2}}\,dx_2''\\
&& + \frac{\cos\theta_3\sqrt{{x_1''}^2+{x_2''}^2}}{a^2}\,dx_3'', \\
\zeta^{(2)}_5 \!\!\!\!&=&\!\!\!\! -\frac{x_1''}{a}\,dx_1'-\frac{x_2''}{a}\,dx_2'-\frac{x_3''}{a}\,dx_3', \\
\zeta^{(2)}_6 \!\!\!\!&=&\!\!\!\!
\frac{x_1'x_1''+x_2'x_2''+x_3'x_3''}{a}\,dt-\frac{x_1''}{a}\,dx_1-\frac{x_2''}{a}\,dx_2-\frac{x_3''}{a}\,dx_3, \\
\zeta^{(2)}_7 \!\!\!\!&=&\!\!\!\!
\frac{ax_2''\sin\theta_3-x_1''x_3''\cos\theta_3}{a\sqrt{{x_1''}^2+{x_2''}^2}}\,dx_1-\frac{ax_1''\sin\theta_3+x_2''x_3''\cos\theta_3}{a\sqrt{{x_1''}^2+{x_2''}^2}}\,dx_2\\
& & -\frac{\sqrt{{x_1''}^2+{x_2''}^2}\cos\theta_3}{a} \,dx_3, \\
\zeta^{(2)}_8 \!\!\!\!&=&\!\!\!\!
\frac{ax_2''\cos\theta_3-x_1''x_3''\sin\theta_3}{a\sqrt{{x_1''}^2+{x_2''}^2}}\,dx_1+\frac{-ax_1''\cos\theta_3+x_2''x_3''\sin\theta_3}{a\sqrt{{x_1''}^2+{x_2''}^2}}\,dx_2\\
& & + \frac{\sqrt{{x_1''}^2+{x_2''}^2}\sin\theta_3}{a} \,dx_3,\\
\zeta^{(2)}_9 \!\!\!\!&=&\!\!\!\! \frac {
a(x_2'x_1''\!-\!x_1'x_2'')\cos\theta_3
\!-\!(x_2'x_3''x_2''\!+\!x_1'x_3''x_1''\!+\!x_3'{x_1''}^{2}\!+\!x_3'{x_2''}^{2})\sin\theta_3
}{a\sqrt{{x_1''}^2+{x_2''}^2}}\,dt \\
& & +
\frac{ax_2''\cos\theta_3+x_1''x_3''\sin\theta_3}{a\sqrt{{x_1''}^2+{x_2''}^2}}\,dx_1
-\frac{ax_1''\cos\theta_3-x_2''x_3''\sin\theta_3}{a\sqrt{{x_1''}^2+{x_2''}^2}}\,dx_2\\
& & +\frac{\sin\theta_3\sqrt{{x_1''}^2+{x_2''}^2}}{a}\,dx_3,\\
\zeta^{(2)}_{10} \!\!\!\!&=&\!\!\!\! \frac {
a(x_2'x_1''\!-\!x_1'x_2'')\sin\theta_3
\!+\!(x_2'x_3''x_2''\!+\!x_1'x_3''x_1''\!+\!x_3'{x_1''}^{2}\!+\!x_3'{x_2''}^{2})\cos\theta_3
}{a\sqrt{{x_1''}^2+{x_2''}^2}}\,dt\\
& & +
\frac{ax_2''\sin\theta_3-x_1''x_3''\cos\theta_3}{a\sqrt{{x_1''}^2+{x_2''}^2}}\,dx_1
-\frac{ax_1''\sin\theta_3+x_2''x_3''\cos\theta_3}{a\sqrt{{x_1''}^2+{x_2''}^2}}\,dx_2\\
& & -\frac{\cos\theta_3\sqrt{{x_1''}^2+{x_2''}^2}}{a}\,dx_3.
\end{eqnarray*}
For any constant $a>0$, these forms serve as a coframe on
\begin{eqnarray*} {\cal B}_2 &=&
\Big\{ (t,{\rm x},{\rm x}_t,{\rm x}_{tt},{\rm
\theta})\in {\Bbb R}^{10} \times {\rm SO}(3) \,\Big|\, \\
& &  \hspace{1cm} \theta_1=-\arctan\Big(\frac{x_2''}{x_1''}\Big),
\theta_2=\arcsin\Big(\frac{x_3''}{\|{\bf x}_{tt}\|}\Big),\|{\rm
x}_{tt}\|=a \Big\}
\end{eqnarray*}}
\section{Third order moving coframes}
By restricting to the third prolongation $J^3{\cal B}_2$ which is
a 2-dimensional manifold, we have,
\paragraph{Theorem 11.}
{\it The restricted (or horizontal) moving coframe forms are
explicitly given by $\eta^{(2)}_1=dt$,
\begin{eqnarray*} \eta^{(2)}_2&=&
\frac{x_3''(x_1'x_2''-x_2'x_1'')}{a({x_1''}^2+{x_2''}^2)}\,\eta^{(2)}_1-d\theta_3,
\end{eqnarray*}
their linear dependencies in this step are
$\eta^{(2)}_3=J_1\eta^{(2)}_1$, $\eta^{(2)}_4=J_2\eta^{(2)}_1$,
$\eta^{(2)}_5=-a\eta^{(2)}_1$ and
$\eta^{(2)}_6=\eta^{(2)}_7=\eta^{(2)}_8=\eta^{(2)}_9=\eta^{(2)}_{10}=0$;
Where
\begin{eqnarray*}
J_1 &=& \frac{1}{a^2\sqrt {{x_1''}^2+{x_2''}^2}}\Big\{
a(x_1^{(3)}x_2''-x_2^{(3)}x_1'')\cos\theta_3 \\
& & + \big((x_1''x_1^{(3)}+x_2^{(3)}x_2'')x_3''
-({x_1''}^{2}+{x_2''}^2)x_3^{(3)}\big)\sin\theta_3
 \Big\}, \\
J_2 &=& \frac{1}{a^2\sqrt {{x_1''}^2+{x_2''}^2}}\Big\{
a(x_1^{(3)}x_2''-x_2^{(3)}x_1'')\sin\theta_3\\
&&
+\big(({x_1''}^{2}+{x_2''}^2)x_3^{(3)}-(x_1''x_1^{(3)}+x_2^{(3)}x_2'')x_3''\big)\cos\theta_3\Big\}.
\end{eqnarray*}}
If we assume $J_1=0$, then we find that $\displaystyle
J_2=\frac{\|{\rm x}_{tt}\times{\rm x}_{ttt}\|}{a^2}$ and
$$\theta_3=\arctan\left(\frac{a(x_1^{(3)}x_2''-x_2^{(3)}x_1'')}{
({x_1''}^{2}+{x_2''}^2)x_3^{(3)}-(x_1''x_1^{(3)}+x_2^{(3)}x_2'')x_3''}\right).$$
thus,
\paragraph{Theorem 12.}
{\it The most general third order compatible lift has the form
\begin{eqnarray*}
\rho^{(3)}(t,{\bf x};{\bf x}_t,{\bf x}_{tt},{\bf x}_{ttt}) =
\left[\begin{array}{ccccc} 1 & 0 & 0 & 0 & t \\ {\rm x}_t &
\frac{{\rm x}_{tt}}{\|{\rm x}_{tt}\|} & \frac{{\rm
x}_{tt}\times{\rm x}_{ttt}}{\|{\rm x}_{tt}\times{\rm x}_{ttt}\|} &
\frac{{\rm x}_{tt}\times({\rm x}_{tt}\times{\rm x}_{ttt})}{\|{\rm
x}_{tt}\times({\rm x}_{tt}\times{\rm x}_{ttt})\|} & {\rm x} \\
0 & 0 & 0 & 0 & 1
\end{array}\right].
\end{eqnarray*}}
\paragraph{Theorem 13.} {\it The third order moving coframe
${\bf \zeta}^{(3)}=\rho^{(3)*}{\bf \mu}$ is
\begin{eqnarray*}
\zeta^{(3)}_1 &=& dt,\\
\zeta^{(3)}_2 &=& \frac{{\rm x}_{tt}\times {\rm x}_{ttt}}{a\|{\rm
x}_{tt}\times {\rm x}_{ttt}\|^2}\cdot\Big(({\rm x}_{tt}\cdot {\rm
x}_{ttt})\,d{\rm x}_{tt} + a^2\,d{\rm x}_{ttt}\Big),\\
\zeta^{(3)}_3 &=& -\frac{{\rm x}_{tt}\cdot {\rm x}_{ttt}}{a\|{\rm
x}_{tt}\times {\rm x}_{ttt}\|}\,({\rm x}_{tt}\times {\rm
x}_{ttt})\cdot d{\rm x}_{tt},\\
\zeta^{(3)}_4 &=& -\frac{{\rm x}_{tt}\times({\rm x}_{tt}\times
{\rm
x}_{ttt})}{a^2\|{\rm x}_{tt}\times {\rm x}_{ttt}\|}\cdot d{\rm x}_{tt},\\
\zeta^{(3)}_5 &=& - \frac{1}{a}{\rm x}_{tt}\cdot d{\rm x}_{t},\\
\zeta^{(3)}_6 &=& \frac{{\rm x}_{t}\cdot {\rm x}_{tt}}{a}\,dt - \frac{1}{a}{\rm x}_{tt}\cdot d{\rm x},\\
\zeta^{(3)}_7 &=& \frac{{\rm x}_{tt}\times({\rm x}_{tt}\times {\rm
x}_{ttt})}{a\|{\rm x}_{tt}\times {\rm x}_{ttt}\|}\cdot d{\rm
x}_{t},\\
\zeta^{(3)}_8 &=& \frac{{\rm x}_{tt}\times {\rm x}_{ttt}}{\|{\rm
x}_{tt}\times {\rm x}_{ttt}\|}\cdot d{\rm x}_{t},\\
\zeta^{(3)}_9 &=& -\frac{{\rm x}_{t}\cdot({\rm x}_{tt}\times {\rm
x}_{ttt})}{\|{\rm x}_{tt}\times {\rm x}_{ttt}\|}\, dt+\frac{{\rm
x}_{tt}\times {\rm x}_{ttt}}{\|{\rm x}_{tt}\times
{\rm x}_{ttt}\|}\cdot d{\rm x},\\
\zeta^{(3)}_{10} &=& \frac{ ({\rm x}_{t}\times {\rm
x}_{tt})\cdot({\rm x}_{tt}\times {\rm x}_{ttt})}{a\|{\rm
x}_{tt}\times {\rm x}_{ttt}\|}\,dt-\frac{{\rm x}_{tt}\times({\rm
x}_{tt}\times {\rm x}_{ttt})}{a\|{\rm x}_{tt}\times {\rm
x}_{ttt}\|}\cdot d{\rm x}.
\end{eqnarray*}
For any constant $a>0$, these forms serve as a coframe on
\begin{eqnarray*}
{\cal B}_3 &=& \Big\{ (t,{\rm x},{\rm x}_t,{\rm x}_{tt},{\rm
\theta})\in {\Bbb R}^{10} \times {\rm SO}(3) \,\Big|\, \\
& & \hspace{1cm} \theta_1=-\arctan\Big(\frac{x_2''}{x_1''}\Big),\;
\theta_2=\arcsin\Big(\frac{x_3''}{\|{\bf x}_{tt}\|}\Big),\|{\rm x}_{tt}\|=a, \\
& & \hspace{1cm}
\theta_3=\arctan\left(\frac{a(x_1^{(3)}x_2''-x_2^{(3)}x_1'')}{
({x_1''}^{2}+{x_2''}^2)x_3^{(3)}-(x_1''x_1^{(3)}+x_2^{(3)}x_2'')x_3''}\right)
\Big\}.
\end{eqnarray*}}
\paragraph{Theorem 14.}
{\it The restricted (or horizontal) moving coframe forms are
explicitly given by $\eta^{(3)}_1=dt$, $\displaystyle
\eta^{(3)}_4=\frac{\|{\rm x}_{tt}\times{\rm
x}_{ttt}\|}{a^2}\,\eta^{(3)}_1$,
$\eta^{(3)}_3=\eta^{(3)}_5=\eta^{(3)}_6=\eta^{(3)}_7=\eta^{(3)}_9=\eta^{(3)}_{10}=0$,
$\eta^{(3)}_8=-a\eta^{(3)}_1$, and $\eta^{(3)}_2=J\eta^{(3)}_1$,
where $J=a\big(({\rm x}_{tt}\times{\rm x}_{ttt})\cdot{\rm
x}_{tttt}\big)/{\|{\rm x}_{tt}\times{\rm x}_{ttt}\|^2}$.}
\paragraph{Theorem 15.} {\it $\displaystyle \frac{d}{dt}$ is a
differential operator and the functions $a_1=\|{\rm x}_{tt}\|$,
$a_2=\|{\rm x}_{tt}\times{\rm x}_{ttt}\|$ and $a_3={({\rm
x}_{tt}\times{\rm x}_{ttt})\cdot{\rm x}_{tttt}}$ are differential
invariants.}
\section{Dimensional considerations}
In this section, we use the conventions of chapter 5 of \cite{O}.

If we use the coordinates $(t,{\rm x},{\rm x}_t,{\rm
x}_{tt},\cdots,{\rm x}^{(n)})$ for $J^n({\Bbb R};{\Bbb R}^3)$,
then the prolonged group action ${\rm SGal}(n)$ on $J^n({\Bbb
R};{\Bbb R}^3)$ can be written as $\bar{t} = t+s$, $\bar{\rm x} =
R\,{\rm x}+t{\rm v}+{\rm y}$, $\bar{\rm x}_{\bar{t}} = R\,{\rm
x}_t+{\rm v}$, and $\bar{\rm x}^{(n)} = R\,{\rm x}^{(n)}$ for
$n\geq2$.

It is recommended that the dimension of $J^n({\Bbb R};{\Bbb R}^3)$
is $p+q^{n}=3n+4$, and the dimension of ${\rm SGal}^{(n)}$ is 10.
\paragraph{Theorem 16.} {\it The following functions are differential
invariants:
\begin{itemize}
\item[1)] $I_n=\|{\rm x}^{(n)}\|$ for $n\geq2$.
\item[2)] $J_{n,m}={\rm x}^{(n)}\cdot{\rm x}^{(m)}$ for $n>m\geq2$.
\item[3)] $K_{n,m}=\|{\rm x}^{(n)}\times{\rm x}^{(m)}\|$ for $n>m\geq2$.
\item[4)] $L_{l,n,m}=({\rm x}^{(l)}\times{\rm x}^{(n)})\cdot{\rm x}^{(m)}$ for $l>n>m\geq2$.
\end{itemize} }

\medskip \noindent {\it Proof:} If $n,m\geq2$, then since $\bar{\rm x}^{(n)} = R\,{\rm x}^{(n)}$, $\bar{\rm x}^{(m)} = R\,{\rm x}^{(m)}$ and $R\in {\rm SO}(3)$,
therefore $\bar{\rm x}^{(n)}\cdot\bar{\rm x}^{(m)}={\rm
x}^{(n)}\cdot{\rm x}^{(m)} $; hence $I_{n,m}$ is an invariant.

By (1), (2) and formulas $\|u\times v\|^2=\|u\|^2\|v\|^2-(u\cdot
v)^2$, we find that $K_{n,m}=I_nI_m-J^2_{n,m}$ is an invariant.

Since $(u_1\times u_2)\cdot u_3 = \det(u_i\cdot u_j)$, then
$L_{l,n,m}$ is a function of $J_{n,m}$'s, and this complete the
proof.\hfill\
$\Box$\\

By usual computations, we find that the maximal dimension of
prolonged action are: $s_0=4$, $s_1=7$, $s_2=9$, $s_n=10$ for
$n\geq3$. Therefore, the order of this group action is $s=3$.

Therefore, the $i_n$ functionally independent differential
invariants of order at most $n$ are: $i_0=i_1=0$, $i_2=1$ and
$i_n=3(n-2)$ for $n\geq3$.

By the theorem 5.31 of \cite{O}, we have
\paragraph{Theorem 17.} {\it The complete system of $3^{\rm rd}$ order differential
invariants of special Galilean group action are $\|{\rm
x}_{tt}\|$, $\|{\rm x}_{ttt}\|$ and ${\rm x}_{tt}\cdot{\rm
x}_{ttt}$. Locally, every 3 $^{rd}$ order differential invariant
of ${\rm SGal}(3)$ can be written as a function of these
differential invariants. \hfill\ $\Box$}
\paragraph{Corollary 1.}
$a_1=I_2,\;\;\; a_2=J_{3,2}$, $a_3={L_{4,3,2}}/{K^2_{3,2}}$.
\hfill\ $\Box$
\paragraph{Theorem 18.} {\it Every differential invariant of
special Galilean group action is a function of $a=\|{\rm
x}_{tt}\|$, $b=\|{\rm x}_{ttt}\|$ and their derivatives with
respect to $t$.}

\medskip \noindent {\it Proof:} According to theorem 17, it is
enough to show that ${\rm x}_{tt}\cdot{\rm x}_{ttt}$ can be
written as a function of $\|{\rm x}_{tt}\|$ and $\|{\rm
x}_{ttt}\|$. But $\displaystyle \frac{1}{2}\frac{d}{dt}\|{\rm
x}_{tt}\|^2={\rm x}_{tt}\cdot{\rm x}_{ttt}$. \hfill\ $\Box$
\section{$\{e\}-$structure}
The necessary condition for local special Galilean equivalence of
two given motions is that the corresponding invariants are the
same. These produce a large amount of necessary conditions.

For sufficient condition of equivalence for coframes ${\it
{\zeta}}$ we can rewrite two-forms $d{\it \zeta_i}$ in terms of
wedge products of the $\it{\zeta_i}$'s. This produces the {\it
structure functions}. There are our {\it original invariants}, by
differentiation from them we have {\it derived invariants}, now we
can construct a large collection of invariants, whose functional
interrelationships provide a necessary condition for equivalence.
It is time to continue by introducing {\it structure invariants}.
This latter structure serves to define the components of the {\it
structure map}. The $s^{th}$ order{\it classifying space} and the
fully regularity condition on $s^{th}$ order {\it structure map}
leads us to the definition of $s^{th}$ order {\it classifying
manifold} ${\cal C}^{(s)}$ due to chapter $8$ in \cite{O}. In view
of the proposition $8.11$  in \cite{O}, necessary conditions for
the (local) equivalence of coframes are that for each $s\ge0$ ,
their $s^{th}$ order classifying manifolds are overlap. Now the
fully regularity conditions provide that these necessary
conditions are also sufficient.

\end{document}